\documentclass[11pt]{article}
\usepackage{amscd, amsmath, amssymb, hyperref}

\title{Divided powers in the Witt ring of symmetric bilinear forms}
\author{Burt Totaro}
\date{  }

\def\Z{\text{\bf Z}}
\def\Q{\text{\bf Q}}
\def\R{\text{\bf R}}

\def\F{\text{\bf F}}
\def\arrow{\rightarrow}

\def\qed{\ QED }

\def\l{\langle \langle}
\def\r{\rangle \rangle}
\def\rank{\text{rank}}

\begin{document}
\maketitle

\newtheorem{theorem}{Theorem}[section]
\newtheorem{corollary}[theorem]{Corollary}
\newtheorem{lemma}[theorem]{Lemma}

Marshall defined divided power operations
in the Witt ring of symmetric bilinear forms over a field
\cite{Marshall}. (In characteristic not 2, it is equivalent
to talk about quadratic forms.) Similar operations have been considered by
Rost, Garibaldi, and Garrel \cite{Rost, Garrel}, \cite[section 19]{Garibaldi}.
On the other hand, it follows from Garibaldi-Merkurjev-Serre's
work on cohomological invariants that all operations on the Witt ring
are essentially 
linear combinations of the exterior power operations
\cite[Theorem 27.16]{GMS}.
In this paper we find the explicit formula for the divided powers
on the Witt ring as a linear combination of exterior powers. The coefficients
are closely related to the ``tangent numbers'', the coefficients
of the power series for $\tan x$, and thus to Bernoulli numbers
\cite[pp.~259--260]{Comtet}. In this way, the existence of divided
powers on the Witt ring of symmetric bilinear forms is explained by an explicit
definition in terms of linear algebra.

The coefficients in the formula
are most important modulo 2, since the Witt
ring $W(k)$ is an $\F_2$-algebra for all fields $k$ in which
$-1$ is a square. We show that the coefficients of divided powers
on the Witt ring in terms of exterior powers simplify modulo 2
to binomial coefficients.
The divided powers on the Witt ring give another construction
of the Rost-Serre-Kahn divided powers on Milnor $K$-theory modulo 2 
\cite{Kahn} via
the Milnor conjecture $I^n/I^{n+1}\cong K_n^M(k)/2$, proved
by Orlov-Vishik-Voevodsky \cite{OVV}. By Vial \cite{Vial},
the divided powers form a basis for 
all operations on the Milnor
$K$-theory of fields.

This work was supported by NSF grant DMS-2054553.
Thanks to Skip Garibaldi and J.-P.~Serre for their comments.
Serre suggested the proof of the first part of Lemma \ref{lambda}.

\section{Notation}
\label{notation}

Let $k$ be a field (any characteristic is allowed).
A symmetric bilinear form over $k$ means a finite-dimensional vector
space with a nondegenerate symmetric bilinear form. For $k$
of characteristic not 2, these can be identified
with quadratic forms.
We write $\langle a_1,\ldots,a_n\rangle$ for the diagonal
form $a_1 x_1y_1+\cdots+a_n x_ny_n$.
The Grothendieck-Witt ring $GW(k)$ is the Grothendieck
group of symmetric bilinear forms over $k$, with addition corresponding
to orthogonal direct sum and multiplication to tensor product.
For $k$ of characteristic not 2,
Witt's cancellation theorem says that two symmetric bilinear forms
over $k$ are isomorphic if and only if they have the same class
in $GW(k)$ \cite[section II.1]{Lam}.
The Witt ring $W(k)$ is the quotient of $GW(k)$ by the subgroup
generated by the hyperbolic plane $\bf H$;
this subgroup is in fact an ideal.
For $k$ of any characteristic,
the bilinear form $\langle 1,-1\rangle$ is zero in $W(k)$,
by the isomorphism $\langle 1,1,-1\rangle\cong \langle 1\rangle\perp
{\bf H}$ \cite[equation 1.16]{EKM}.

We can identify the ideal $I=\ker(W(k)\arrow \Z/2)$
of even-dimensional forms
with $\ker(\rank : GW(k)\arrow \Z)$.
The Grothendieck-Witt ring $GW(k)$ has the advantage
of being a $\lambda$-ring, using
exterior powers of symmetric bilinear forms.
(That $GW(k)$ is a $\lambda$-ring
follows from its being additively generated by 1-dimensional forms
\cite[Proposition 9.8]{McGarraghy}; this argument works
in any characteristic.)
A reference on $\lambda$-rings is \cite{AT}.
As is now standard, we call
a $\lambda$-ring what \cite{AT} calls a special $\lambda$-ring.

Following \cite[Definition 3.1]{BObook},
a divided power
structure on an ideal $I$ in a commutative ring $R$
is a collection of functions $\gamma_n:I\arrow R$ for
$n\geq 0$ such that:

(1) $\gamma_0(x)=1$, $\gamma_1(x)=x$, $\gamma_n(x)\in I$ for $n>0$
and $x\in I$.

(2) $\gamma_n(x+y)=\sum_{i=0}^n \gamma_i(x)\gamma_{n-i}(y)$ for $x,y\in I$.

(3) $\gamma_n(ax)=a^n\gamma_n(x)$ for $a\in R$, $x\in I$.

(4) $\gamma_m(x)\gamma_n(x)=\binom{m+n}{m} \gamma_{m+n}(x)$ for $x\in I$.

(5) $\gamma_n(\gamma_m(x))=\frac{(mn)!}{(m!)^n n!}\gamma_{mn}(x)$ for $x\in I$.
(The coefficient is an integer.)

Here relation (4) implies that $n!\, \gamma_n(x)=x^n$, and so
every ideal in a commutative $\Q$-algebra has a unique divided power
structure defined by $\gamma_n(x)=x^n/n!$; this explains
where the axioms come from.

\section{The formula for divided powers}

We give a formula for the divided power structure on the ideal $I$
of even-dimensional forms
in the Witt ring localized at 2. In general,
it is necessary to localize the Witt ring
at 2 to have divided powers, as one sees by looking at
the Witt ring of the real numbers, $W(\R)=\Z$. Note that localizing at 2
makes very little difference for Witt rings.
Indeed, Pfister showed that
there is no odd torsion in $W(k)$,
and so 
$W(k)$ always injects into $W(k)_{(2)}$. Moreover, if $-1$
is a sum of squares in $k$ (which holds in all fields
of positive characteristic, for example), then $W(k)$ is killed
by some power of 2 and hence is equal to its localization at 2
\cite[Theorem VIII.3.2]{Lam}.

\begin{theorem}
\label{main}
Let $k$ be a field. Then the ideal
$I_{(2)} =\ker(W(k)\arrow \Z/2)_{(2)}$ in the Witt ring localized
at 2 has a divided power structure defined by
$$\gamma_n(x)=\sum_{i\geq 0} (-1)^{(n-i)/2} (i!/n!) T(n,i)\lambda^i(x),$$
where $T(n,i)$ are the ``tangent numbers'' defined by
$$\frac{(\tan t)^i}{i!}=\sum_{n\geq i}T(n,i)t^n/n! .$$
Here we identify $I_{(2)}$ with $\ker(GW(k)\arrow \Z)_{(2)}$, in which
exterior powers are defined.
The coefficients in the formula for $\gamma_n(x)$
are 2-local integers, so that
the formula makes sense in $W(k)_{(2)}$.
\end{theorem}

Note that the divided power operations are not the
``gamma operations'' which exist in any $\lambda$-ring.
The sign in the formula makes sense because $T(n,i)$
is nonzero only for $i\equiv n\pmod{2}$. Various formulas
for tangent numbers are discussed by Comtet \cite[pp.~259--260]{Comtet},
from which the following table is taken. Many references
on tangent numbers can be found in \cite{Sloane}.

$$\begin{tabular}{r|rrrrrrrrrr}
$n\backslash i$ & 0&1&2&3&4&5&6&7&8&9 \\
\hline
0&1 & &&&&&&&& \\
1& & 1 &&&&&&&& \\
2&& & 1 &&&&&&& \\
3&& 2 & & 1 &&&&&& \\
4&& & 8 & & 1 &&&&& \\
5&& 16 & & 20 && 1 &&&& \\
6&&& 136 && 40 && 1 &&& \\
7&& 272 && 616 && 70 && 1 && \\
8&&& 3968 && 2016 && 112 && 1 & \\
9&& 7936 && 28160 && 5376 && 168 && 1
\end{tabular}$$

The numbers in column 1, giving the coefficients
of the Taylor series of the tangent function, are closely related
to Bernoulli
numbers, in the sense that 
$$T(2n+1,1)=(-1)^{n-1}B_{2n}4^n(4^n-1)/2n.$$
Therefore, one cannot expect too simple a formula for the
coefficients of $\gamma_n$ in terms of exterior powers,
for an arbitrary field. The first few formulas are:
\begin{align*}
\gamma_1(x)&=x\\
\gamma_2(x)&=\lambda^2(x),\\
\gamma_3(x)&=\lambda^3(x)-(1/3)x\\
\gamma_4(x)&=\lambda^4(x)-(2/3)\lambda^2(x).
\end{align*}
On the other hand,
for fields in which $-1$ is a square, Corollary \ref{binom} below
gives a simple formula for the divided powers $\gamma_n$.

{\bf Proof of Theorem \ref{main}. }The axioms 
imply that a divided power structure on an ideal $I$
in a commutative $\Z_{(2)}$-algebra $R$ is uniquely determined by
the operation $\gamma_2$. Moreover, a function
$\gamma_2:I\arrow I$ gives a divided power structure exactly
when it satisfies $\gamma_2(x+y)=\gamma_2(x)+xy+\gamma_2(y)$
and $\gamma_2(ax)=a^2\gamma_2(x)$ for $x,y\in I$, $a\in R$.
See Berthelot-Ogus \cite[Appendix]{BOpaper} for the analogous description
of divided power structures in $\Z_{(p)}$-algebras for any prime
number $p$.
(They also assume that $2\gamma_2(x)=x^2$, but that follows
from the formulas mentioned.)

Thus we can define a divided power structure on the ideal
$I_{(2)}=\ker(GW(k)\arrow \Z)_{(2)}$ in $GW(k)_{(2)}$ by declaring
that $\gamma_2(x)=\lambda^2(x)$ for $x\in I_{(2)}$. Here the first identity
$\lambda^2(x+y)=\lambda^2(x)+xy+\lambda^2(y)$ is a standard
property of exterior powers (it holds in any $\lambda$-ring).
The second property is special to the Grothendieck-Witt ring.
Indeed, for all elements $a,x$ of a $\lambda$-ring, we have
$$\lambda^2(ax)=a^2\lambda^2(x)+\lambda^2(a)\psi^2(x),$$
where $\psi^2$ is the Adams operation defined by
$\psi^2(x)=x^2-2\lambda^2(x)$. So it suffices to show
that $\psi^2$ is zero on the ideal $I$ in $GW(k)$. 
Since Adams operations are ring homomorphisms, this is easy.
View any element of $I$ as the difference of two symmetric bilinear
forms of the same dimension, and we have:
\begin{align*}
\psi^2(\langle a_1,\ldots,a_n\rangle -\langle b_1,\ldots,b_n\rangle)&=
\langle a_1^2,\ldots,a_n^2\rangle -\langle b_1^2\ldots,b_n^2\rangle\\
&= \langle 1,\ldots, 1\rangle - \langle 1,\ldots, 1\rangle \\
&= 0.
\end{align*}
Thus we have constructed a divided power structure on
the ideal $I_{(2)}=\ker(GW(k)\arrow \Z)_{(2)}$. This is the same
divided power structure as that constructed by Marshall \cite{Marshall}.
It is immediate that these operations give
a divided power structure on $I_{(2)}$
as an ideal in the Witt ring, $I_{(2)}=\ker(W(k)\arrow \Z/2)_{(2)}$.

From the axioms for a divided power structure, all the operations
$\gamma_n$ on $I_{(2)}$ are $\Z_{(2)}$-polynomials in iterates of $\gamma_2
=\lambda^2$. But the Grothendieck-Witt ring of a field
has extra properties
not valid in an arbitrary $\lambda$-ring, using that it is generated
by elements $\alpha$ with $\alpha^2=1$ and $\lambda^i(\alpha)=0$
for $i>1$.  In particular, we show in the proof
of Lemma \ref{lambda} that
$\lambda^a(x)\lambda^b(x)$ is a $\Z$-linear combination
of the exterior powers $\lambda^m(x)$ for $m\leq a+b$,
and $\lambda^a(\lambda^b(x))$ is a $\Z$-linear combination
of $\lambda^m(x)$ for $m\leq ab$ (with coefficients depending
on the rank of $x$).
Therefore, for $x$ of rank zero, we have
$$\gamma_n(x)=\sum_{i\leq n}a(n,i)\lambda^i(x)$$
for some universal coefficients $a(n,i)$ in $\Z_{(2)}$
(independent of the field $k$). We want to determine these coefficients.

We first give the formulas for $\lambda^a(x)\lambda^b(x)$
and $\lambda^2(\lambda^a(x))$.
(In this paper, we only need the case
of elements $x\in GW(k)$ of rank zero.)

\begin{lemma}
\label{lambda}
Let $k$ be a field.
For an element $x\in GW(k)$ of rank $N\in\Z$,
$$\lambda^a(x)\lambda^b(x)=\sum_{0\leq j\leq \min(a,b)} \binom{a+b-2j}{a-j}
\binom{N-a-b+2j}{j}\lambda^{a+b-2j}(x)$$
and
$$\lambda^2(\lambda^a(x))=\sum_{0\leq j<a}\frac{1}{2}
\binom{2a-2j}{a-j}\binom{N-2a+2j}{j} \lambda^{2a-2j}(x).$$
The coefficients in these formulas are integers.
\end{lemma}

{\bf Proof of Lemma \ref{lambda}.}
For an element $x\in GW(k)$, write
$$\lambda_t(x)=1+tx+t^2\lambda^2(x)+\cdots$$
in the power series ring $GW(k)[[t]]$. Let $N\in \Z$ be the rank
of $x$. Serre showed that
$$\lambda_u(x)\lambda_v(x)=(1+uv)^N\lambda_{(u+v)/(1+uv)}(x)$$
in $GW(k)[[u,v]]$ \cite[Exercise 27.2(3)]{GMS}. (Use that $\lambda_t(x+y)
=\lambda_t(x)\lambda_t(y)$ to reduce to the case where $x$ is a 1-dimensional
form, in which case the formula is clear.) Equivalently,
$$\sum_{a,b\geq 0}\lambda^a(x)\lambda^b(x)u^av^b
=\sum_{m\geq 0}(u+v)^m(1+uv)^{N-m}\lambda^m(x).$$
We have $(u+v)^m$ as $\sum_{0\leq c\leq m}\binom{m}{c}u^cv^{m-c}$,
and likewise $(1+uv)^{N-m}=\sum_{j\geq 0}\binom{N-m}{j}u^jv^j$. So
$$\sum_{a,b\geq 0}\lambda^a(x)\lambda^b(x)u^av^b=\sum_{c,m,j\geq 0}
\binom{m}{c}\binom{N-m}{j}u^{c+j}v^{m-c+j}\lambda^m(x).$$
On the right side, the $u^av^b$ terms are obtained when
$c=a-j$ and $m=a+b-2j$.
That gives the desired
formula for $\lambda^a(x)\lambda^b(x)$.

Since $GW(k)$ is a $\lambda$-ring,
$\lambda^a(\lambda^b(x))$
is given by a universal polynomial with integer coefficients
\cite[Definition 9.6]{McGarraghy}:
$$\lambda^a(\lambda^b(x))=Q_{a,b}(\lambda^1x,\ldots,\lambda^{ab}x).$$
Explicitly, let $\xi_1,\ldots,\xi_n$ be variables, with $n\geq ab$.
Let $e_1,e_2,\ldots$
be the elementary symmetric functions in $\xi_1,\ldots,\xi_n$.
Then
$$Q_{a,b}(e_1,\ldots,e_{ab}):=\text{coefficient of }t^a\text{ in }
\prod_{1\leq i_1<\cdots<i_b\leq n}(1+\xi_{i_1}\cdots \xi_{i_b}t).$$
Clearly the polynomial $Q_{a,b}$ is homogeneous of degree $ab$,
with each $e_i$ viewed as having degree $i$. Therefore, in all
Grothendieck-Witt rings of fields, the formula above
for $\lambda^a(x)\lambda^b(x)$ gives a universal formula
$$\lambda^a(\lambda^b(x))=\sum_{j\geq 0}f_{a,b,j}(N)\lambda^{ab-2j}(x),$$
with each $f_{a,b,j}$ an integer-valued polynomial in $N:=\rank(x)$.

To determine the coefficients in one of these formulas, we can work in a
field $k$ in which $1,x,\lambda^2(x),\ldots,\lambda^r(x)$
(for a given number $r$) are linearly
independent in $GW(k)\otimes\Q$ for some $x\in GW(k)$ (or, equivalently,
$1,x,x^2,\ldots,x^r$ are linearly independent in $GW(k)\otimes\Q$).
For example, this holds
for $k$ a suitable rational function field over the real numbers.
In this situation, we do not lose any information
by tensoring $GW(k)$ with the rationals.

To work out the formula for $\lambda^2(\lambda^a(x))$, we use
again that the Adams operation $\psi^2(x)=x^2-2\lambda^2(x)$
is a ring homomorphism. For a 1-dimensional form $\langle c\rangle$, we have
$\psi^2\langle c\rangle=\langle c^2\rangle =1$; so, for any $x\in GW(k)$,
$\psi^2(x)$ is equal to the rank $N$ of $x$. That is,
$\lambda^2(x)=(-N+x^2)/2$. Since $\lambda^a(x)$ has rank $\binom{N}{a}$,
we can compute $\lambda^2(\lambda^a(x))$ from the formula above
for $(\lambda^a(x))^2$:
\begin{align*}
\lambda^2(\lambda^a(x))&=\frac{1}{2}\bigg[-\binom{N}{a}
+\sum_{0\leq j\leq a}\binom{2a-2j}{a-j}\binom{N-2a+2j}{j}\lambda^{2a-2j}(x)\bigg]\\
&=\sum_{0\leq j<a}\frac{1}{2}\binom{2a-2j}{a-j}\binom{N-2a+2j}{j}
\lambda^{2a-2j}(x),
\end{align*}
as we want. \qed

\medskip

Returning to the proof of Theorem \ref{main},
let $x$ be an element of $I=\ker(GW(k)\to\Z)$. We know that there
is a formula for the divided power $\gamma_n(x)$
as a $\Z_{(2)}$-linear combination
of the elements $\lambda^a(x)$.
As in the proof of Lemma \ref{lambda},
we can work in a
field $k$ in which $1,x,\lambda^2(x),\ldots,\lambda^n(x)$ are linearly
independent in $GW(k)\otimes\Q$ for some $x\in I$ (or, equivalently,
$1,x,x^2,\ldots,x^n$ are linearly independent in $GW(k)\otimes\Q$).
In this situation, we do not lose any information
by tensoring $GW(k)$ with the rationals. Thus
the divided powers are simply
$\gamma_n(x)=x^n/n!$. In particular, we have
$\gamma_{n+1}(x)=x\gamma_n(x)/(n+1)$
for all $n\geq 0$. Also, we have
$$x\lambda^i(x)=(i+1)\lambda^{i+1}(x)-(i-1)\lambda^{i-1}(x)$$
by Lemma \ref{lambda}. This gives
a recurrence relation for the numbers $a(n,i)\in \Z_{(2)}$:
$$a(n+1,i)=\frac{i}{n+1}(a(n,i-1)-a(n,i+1)).$$

Since $\gamma_0(x)=1$, the numbers $a(0,i)$ in the 0th row
are 1 for $i=0$
and $0$ otherwise. By the recurrence relation, we have
$a(n,i)=0$ unless $n\equiv i\pmod{2}$. The
statement we are trying to prove suggests defining
rational numbers $b(n,i)$ by
$a(n,i)=(-1)^{(n-i)/2} (i!/n!) b(n,i)$ for $i,n\geq 0$.
The recurrence relation for $a(n,i)$ implies that
$$b(n+1,i)=b(n,i-1)+i(i+1)b(n,i+1).$$
But this is exactly the recurrence relation satisfied by the
tangent numbers $T(n,i)$ defined by
$$(\tan t)^i/i! = \sum_{n\geq i}T(n,i)t^n/n!$$
\cite[p.~259]{Comtet}. To check that, differentiate
this formula for $(\tan t)^i/i!$, using that
the derivative of $\tan t$ is $1+(\tan t)^2$. \qed (Theorem \ref{main})

\section{Divided powers when $-1$ is a square}

We now show that the coefficients in the formula
for divided powers in the Witt ring simplify to binomial
coefficients modulo 2. This is relevant to fields $k$
in which $-1$ is a square, since then $W(k)$
is an $\F_2$-algebra.

\begin{corollary}
\label{binom}
Let $k$ be a field in which $-1$ is a square. Then
the ideal $I=\ker(W(k)\arrow \Z/2)$ has a divided power
structure defined by the formula
$$\gamma_n(x)=\sum_j \binom{n}{j}\lambda^{n-2j}(x).$$
Here we identify $I$ with $\ker(GW(k)\arrow \Z)$, in which
exterior powers are defined.
\end{corollary}

For example, when $-1$ is a square in $k$, we have:
\begin{align*}
\gamma_1(x)&=x\\
\gamma_2(x)&=\lambda^2(x)\\
\gamma_3(x)&=\lambda^3(x)+x\\
\gamma_4(x)&=\lambda^4(x).
\end{align*}

{\bf Proof. }As shown in the proof of Theorem \ref{main},
$I_{(2)}=\ker(GW(k)\arrow \Z)_{(2)}$
has a unique divided power structure such that
$\gamma_2(x)=\lambda^2(x)$. Since we now assume
that $-1$ is a square in $k$, $I=I_{(2)}$ is killed by 2.
By repeated application of the formula for $\gamma_n(\gamma_m(x))$
in a divided power ideal, it follows that
$\gamma_{2^r}(x)=(\gamma_2)^r(x)$ for $x\in I$.

Next, by Lemma \ref{lambda},
$$\lambda^2(\lambda^{2^{r}}(x))=\sum_{j\geq 0}\frac{1}{2}
\binom{2^{r+1}-2j}{2^{r}-j}\binom{-2^{r+1}+2j}{j} \lambda^{2^{r+1}-2j}(x).$$
Consider the first binomial coefficient here:
$(1/2)\binom{2n}{n}$ is 0 modulo 2 except when $n$ is a power of 2,
where it is 1 modulo 2.
Since we are working modulo 2 (as $-1$ is a square in $k$), most terms
in the sum disappear and we have
$$\lambda^2(\lambda^{2^{r}}(x))=\sum_{s=0}^r \binom{-2^{s+1}}
{2^r-2^s}\lambda^{2^{s+1}}(x).$$
The coefficient here is 0 modulo 2 for $s<r$ and 1 modulo 2 for $s=r$.
We conclude that
$$\lambda^2(\lambda^{2^r}(x))=\lambda^{2^{r+1}}(x).$$
Therefore, by induction,
$\gamma_{2^r}(x)=(\lambda^2)^r(x)=\lambda^{2^r}(x)$ for $x\in I$. This proves
Corollary \ref{binom} for $n=2^r$, since $\binom{2^r}{j}=0\pmod {2}$
for $0<j<2^r$.

Let $n=2^{i_0}+\cdots+2^{i_r}$ be the binary expansion of $n$.
Then $\gamma_n(x)$ is a constant in $\Z_{(2)}^*$ times
$\gamma_{2^{i_0}}(x)\cdots\gamma_{2^{i_r}}(x)$ for $x\in I$. Since
we can work modulo 2,
\begin{align*}
\gamma_n(x)&=\gamma_{2^{i_0}}(x)\cdots\gamma_{2^{i_r}}(x)\\
&= \lambda^{2^{i_0}}(x)\cdots\lambda^{2^{i_r}}(x).
\end{align*}
We want to show that this equals $\sum_{j\geq 0}
\binom{n}{j}\lambda^{n-2j}(x)$. 

We prove this formula by induction
on the number of ones in the binary expansion of $n$. Thus,
we suppose that the formula holds for $m:=2^{i_1}+\cdots + 2^{i_r}$,
and we will prove it for $n=2^{i_0}+m$, where $i_0<i_1< \cdots
< i_r$. We have, for $x\in I$,
\begin{align*}
\gamma_n(x) &= \lambda^{2^{i_0}}(x)\cdots\lambda^{2^{i_r}}(x) \\
&= \lambda^{2^{i_0}}(x)\sum_{j\geq 0}\binom{m}{j}\lambda^{m-2j}(x)\\
&= \sum_{j,l\geq 0} \binom{m}{j} \binom{2^{i_0}+m-2j-2l}{2^{i_0}-l}
\binom{-2^{i_0}-m+2j+2l}{l}\lambda^{n-2j-2l}(x),
\end{align*}
using Lemma \ref{lambda}.
Here $\binom{m}{j}$ is 0 modulo 2 unless $2^{i_1}|j$,
since $m$ is a multiple of $2^{i_1}$; so we can assume that
$2^{i_1}$ divides $j$ in the sum. So $2^{i_0}|(2^{i_0}+m-2j)$.
If $2^{i_0}\nmid l$, then the bottom number in the binomial
coefficient $\binom{2^{i_0}+m-2j-2l}{2^{i_0}-l}$ has lowest
binary digit 1 where the top number has digit 0, and
so the binomial coefficient is 0 modulo 2. So we can also assume
that $2^{i_0}| l$ in the above sum. But the binomial
coefficient just mentioned is also 0 if $l>2^{i_0}$, and so
we have $l=0$ or $l=2^{i_0}$.  That is,
$$\gamma_n(x)=\sum_{j\geq 0} \binom{m}{j}\binom{n-2j}{2^{i_0}}
\lambda^{n-2j}(x)+\sum_{j\geq 0}\binom{m}{j} \binom{2^{i_0}-m+2j}
{2^{i_0}}\lambda^{n-2^{i_0+1}-2j}(x).$$

As we mentioned, the terms in these sums are 0 unless
$2^{i_1}|j$. Given that $2^{i_1}|j$, we have
$n-2j\equiv 2^{i_0}\pmod{2^{i_0+1}}$, and hence the binomial
coefficient
$\binom{n-2j}{2^{i_0}}$ in the left sum is 1 modulo 2
for $0\leq j\leq n/2$
and $2^{i_1}|j$.  (We only need to consider $j$ at most $n/2$
since we are studying the coefficient of $\lambda^{n-2j}(x)$.)
Likewise for the negative binomial coefficient in the right sum, above:
$2^{i_1}$ divides $m-2j$, and so $\binom{2^{i_0}-m+2j}{2^{i_0}}$ 
is 1 modulo 2 for $0\leq j\leq n/2$ and $2^{i_1}|j$. Thus
\begin{align*}
\gamma_n(x)&=\sum_{j\geq 0}\binom{m}{j}\lambda^{n-2j}(x)
+\sum_{j\geq 0}\binom{m}{j}\lambda^{n-2^{i_0+1}-2j}(x)\\
&=\sum_{j\geq 0}\bigg[ \binom{m}{j}+\binom{m}{j-2^{i_0}}
\bigg] \lambda^{n-2j}(x).
\end{align*}

The coefficient in the last sum
is 0 modulo 2 unless $2^{i_0}|j$; note that the index $j$ in this sum
is either the $j$ in the previous sums, which is a multiple
of $2^{i_1}$, or else the old $j$ plus $2^{i_0}$. We know
that $\binom{u}{v}+\binom{u}{v-1}=\binom{u+1}{v}$. Since
multiplying the top and bottom numbers by a power of 2
does not change a binomial coefficient modulo 2, we
have $\binom{2^{i}u}{2^{i}v}+\binom{2^{i}u}
{2^{i}v-2^{i}}\equiv \binom{2^{i}u+2^{i}}{2^{i}v}\pmod{2}$.
Since $m+2^{i_0}=n$, this simplifies our formula for
$\gamma_n(x)$ to:
$$\gamma_n(x)=\sum_{j\geq 0}\binom{n}{j}\lambda^{n-2j}(x).$$
This completes the inductive proof of this formula
when $-1$ is a square in $k$. \qed

\section{Comparison with Milnor $K$-theory}

It is straightforward to compute the divided powers of
Pfister forms in the Witt ring. (By definition, a 1-fold
Pfister form $\l a \r$, for $a\in k^*$, means the 2-dimensional
symmetric bilinear form $\langle 1,-a\rangle$, and an $r$-fold Pfister
form $\l a_1,\ldots,a_r\r$ means the tensor product
$\l a_1\r\cdots \l a_r\r$.) For example, one checks
by induction on $r$ that
\begin{align*}\gamma_2\l a_1,\ldots,a_r\r &=\l a_1,\ldots,a_r\r 
  \l -1\r ^{r-1}\\
&= 2^{r-1}\l a_1,\ldots,a_r\r .
\end{align*}
In particular, if $-1$ is a square
in $k$ (so that $2=0$ in $W(k)$), then
$$\gamma_2\l a_1,\ldots,a_r\r=\begin{cases} \l a_1 \r & \text{if }
r=1\\
0 & \text{if }r\geq 2.
\end{cases}$$
As a result, when $-1$ is a square, we get simple formulas
for the divided powers of Pfister forms: for $n>0$,
$$\gamma_n\l a \r=\begin{cases} \l a \r & \text{if }n\text{ is a power of 2}\\
0 & \text{otherwise}
\end{cases}$$
and, for $r\geq 2$ and $n>0$,
$$\gamma_n\l a_1,\ldots,a_r \r=\begin{cases} \l a_1,\ldots,a_r
 \r & \text{if }n=1\\
0 & \text{otherwise}
\end{cases}$$

This calculation plus the formal properties of divided powers imply
a simple formula for the divided powers of any element
of $I^r$, written as a sum of $r$-fold Pfister forms $s_i$,
when $r\geq 2$ and $-1$ is a square:
$$\gamma_n\bigg( \sum_{i=1}^r s_i\bigg) = \sum_{1\leq i_1<\cdots <i_n\leq r}
  s_{i_1}\cdots s_{i_n}.$$
This implies:

\begin{corollary}
If $-1$ is a square in a field $k$, then the divided power
$\gamma_n$ maps $I^r\subset W(k)$ into $I^{nr}$ for $r\geq 2$.
This operation is compatible with the
divided power operation on Milnor $K$-theory modulo 2,
 $\gamma_n: K^M_r(k)/2=I^r/I^{r+1}\arrow K^M_{nr}(k)/2=I^{nr}/I^{nr+1}$.
\end{corollary}

Indeed, divided powers on Milnor $K$-theory modulo 2
are defined exactly when $-1$ is a square and $r\geq 2$,
and in that case they are defined by the same formula
as above (where $s_i$ are symbols $\{ a_1,\ldots,a_r\}$ in $K^M_r(k)/2$)
\cite{Vial}. (In Kahn's construction of mod 2 divided powers
\cite[Theorem 2]{Kahn}, note that divided powers are not defined
on the whole exterior algebra in characteristic 2, but only in degrees
at least 2.)


\small \sc UCLA Mathematics Department, Box 951555,
Los Angeles, CA 90095-1555

totaro@math.ucla.edu
\end{document}